\documentclass{ifacconf}

\usepackage{graphicx}   
\usepackage{natbib}    
\usepackage{amsmath} 
\usepackage{amssymb} 
\usepackage{color} 
\usepackage{url}

\usepackage{cleveref}

\crefname{equation}{}{} 

\def\reft{{{\textup{ref}}}}

\usepackage{soul} 
\usepackage{xcolor}

\def\mbR{\mathbb{R}}

\def\mcC{\mathcal{C}}

\newcommand{\lonew}{{$\mathcal{L}_1$}}

\DeclareMathOperator*{\trace}{tr}

\def\lone{${\mathcal{L}_1}$ }
\newcommand{\II}{\mathbb{I}}
\newcommand{\IR}{\mathbb{R}}
\newcommand{\norm}[1]{\left\lVert #1 \right\rVert}
\newcommand{\Linfnorm}[1]{ \| #1 \|_{\mathcal{L}_\infty}}
\newcommand{\Proj}[2]{\text{Proj}\left({#1},{#2}\right)}

\newcommand{\abs}[1]{ | #1 |}

\begin{document}
	\begin{frontmatter}
		
		\title{\hspace{-4mm} \lone Adaptive Control for Switching Reference Systems: Application to Flight Control\thanksref{footnoteinfo}\thanksref{footnoteinfo2}} 
		
		\thanks[footnoteinfo]{This work was supported by NASA Langley Research Center.}
		\thanks[footnoteinfo2]{This paper is an updated version of \cite{snyder2019switchingl1} with revised stability condition in Lemma~4 and corrected proof for Theorem 9.}
		
		\author[First]{Steven Snyder} 
		\author[Second]{Pan Zhao} 
		\author[Third]{Naira Hovakimyan} 
		
		\address{University of Illinois Urbana-Champaign, 
			Champaign, IL 61820 USA}
		\addtocounter{address}{-1}
		\address[First]{(e-mail: ssnyder6@illinois.edu)}
		\address[Second]{(e-mail: panzhao2@illinois.edu)}
		\address[Third]{ (e-mail: nhovakim@illinois.edu)}

		\begin{abstract}        
This paper presents a framework for the design and analysis of an $\mathcal{L}_1$ adaptive controller with a switching reference system. The use of a switching reference system allows the desired behavior to be scheduled across the operating envelope, which is often required in aerospace applications. The analysis uses a switched reference system that assumes perfect knowledge of uncertainties and uses a corresponding non-adaptive controller. Provided that this switched reference system is stable, it is shown that the closed-loop system with unknown parameters and disturbances and the $\mathcal{L}_1$ adaptive controller can behave arbitrarily close to this reference system. Simulations of the short period dynamics of a transport class aircraft during the approach phase illustrate the theoretical results.
		\end{abstract}
		
		\begin{keyword}
			Adaptive control, Control System Analysis, Flight Control, Switching System.
		\end{keyword}
		
	\end{frontmatter}
	\section{Introduction}
	In many aerospace applications, a local linear approximation of the plant is sufficient for local control. 
	However, over the entire flight envelope, the dynamics can differ significantly.
	For example, at high speeds, with more air flowing over the control surfaces, there can be significantly more control authority than at low speeds.  
	More modification of the natural dynamics of the vehicle is then possible, so the control objectives can be more ambitious. 
	Sometimes the overall behavior of the system can vary throughout the envelope from stable to unstable. 
	Typically, once control laws are designed in the different local regions of the flight envelope, they are then scheduled throughout. 
	There are many examples of gain scheduling.
	The design for the high-alpha research vehicle by \cite{harv} uses gain scheduling. 
	In \cite{busjet} gain scheduling is used for control of a business jet, and the F-35 control law uses a dynamic inversion controller based on scheduled linear models, per \cite{F35}.
	Applying gain scheduling enables the use of analysis techniques and design criteria intended for linear plants while allowing different dynamics to be set at different locations within the flight envelope. 
	
	In recent years, researchers have begun investigating adaptive control designs for piecewise linear (gain scheduled) systems. 
	In \cite{Sang2012}, piecewise linear reference models are used in the model reference adaptive control (MRAC) framework with projection based adaptation laws. The stability condition is given based on a dwell time argument. 
	The authors of \cite{Yuan2016} sought to extend the results of \cite{Sang2012} by proposing new stability criteria allowing for the Lyapunov matrix to be time varying. However, the adaptive laws and switching laws are coupled in this approach. Both of these results are only applicable to single-input systems. The work in \cite{Kersting2017} provides the first results for multi-input MRAC for piecewise affine systems. 
	
	In the past decade, \lone adaptive control has been developed (see \cite{L1Book}) and applied to numerous flight control applications, including NASA's AirSTAR in \cite{AirSTAR} and Calspan's variable-stability Learjet in \cite{Learjet2019}.
	However, in these flight applications, the control law was either not scheduled or the scheduling was done in an \emph{ad-hoc} fashion, lacking a rigorous mathematical analysis and relying on extensive numerical simulations for stability and performance verification.
	This paper provides a method for analyzing an \lone adaptive controller where the desired dynamics are changing throughout the flight envelope. 
	Note that these changes could occur due to known scheduling parameters, e.g. fuel state, or due to online model identification.

	\section{Motivating Examples} \label{sec::motivation}
	Before presenting the main results, we present a few motivating examples.
	\lone adaptive control has been shown to compensate for uncertainties and disturbances quite well (see e.g. \cite{Kasey2016,AirSTAR}).
	However, sometimes changes in the system response are also important. 
	
	\emph{Changes in Airspeed} \nopagebreak
	
	In \cite{Kasey2019}, during flight testing of an \lone adaptive control law, pilots noted that they could not feel the typical change in stick force associated with the vehicle slowing down, causing pilots to spend more time looking at the instrument gauges.
	As noted by the authors, ``the lack of cuing through stick sensitivity is a direct (expected) result of the adaptive augmentation providing automatic compensation for the change in stick sensitivity and providing a consistent aircraft response despite the deviation in airspeed from the design condition.''
	
	\emph{Changes in Inertia} \nopagebreak
	
	Another common scenario where dynamics are expected to change is when the inertia changes, such as when payload is dropped or fuel stored on the wings of a vehicle is burned. 
	For a given speed, \cite{Learjet2019} shows that the roll mode time constant changes by a factor of roughly 2.5 as fuel is burned.
	If the natural dynamics provide an acceptable response, it may not be worth spending the vehicle's finite control power to modify the dynamics back to a fixed reference system. 
	This would reduce the available control power for disturbance rejection and command tracking while making it more challenging to achieve desired stability margins.

	\emph{System Identification} \nopagebreak
	
	In \cite{L2F}, a new paradigm in aircraft design and testing is suggested. 
	The central idea is to use state-of-the-art system identification techniques 
	to develop mathematical models of the aircraft onboard and in real time during flight. 
	The control algorithms in this framework are adjusted based on the identified models.
	
	The present work aims to add to the \lone adaptive control literature an approach for handling switched reference systems in order to achieve expected (desired) changes in the system dynamics.  
	
	Throughout this paper, 
	we use $\norm{\cdot}$ to denote either the Euclidean norm of a vector or the induced 2-norm of a matrix. 
	$\IR^n$ denotes the $n$-dimensional real vector space. 
	$\II$ denotes an identity matrix of appropriate dimensions. 
	For symmetric matrices $P$ and $Q$, $P>Q$ means $P-Q$ is positive definite. 
	For a function of time $x \colon t\rightarrow \mathbb{R}^n$, its Laplace transform is denoted by $x(s)\triangleq\mathfrak{L}[x(t)]$, and its $\mathcal{L}_\infty$~norm is defined as $\Linfnorm{x}\triangleq \max_{t\geq 0} \norm{x(t)}_2$. 
	
	\section{\lone Adaptive Control for Switching Reference Systems}
	\subsection{Problem Formulation}
	
	Consider the family of multi-input multi-output LTI systems whose state-space matrices are given by 
		\begin{align}
			\{ (A_i, B_i, C_i): i \in \mathcal{I} \}, \label{eq:sub-system}
		\end{align}
		where $\mathcal{I}$ denotes the index set. Let $\mathcal{P} = \{ p : [0,\infty) \to \mathcal{I} \}$ denote the family of piecewise constant switching signals. For a given switching signal $p \in \mathcal{P}$, define the following switched linear system subjected to time-varying parametric uncertainty and disturbances:
	\begin{align}
		\begin{split} \label{eq::ch4sysdyn2}
			\dot{x}(t) &= A_{p} x(t) + B_p\left(\omega_p u(t) + \theta^\top_p(t)x(t) + \sigma_p(t) \right), \\
			y(t) &= C_p x(t), \quad x(0) = x_0, 
		\end{split}
	\end{align}
	where $x(t) \in \IR^n$ is the system state, $u(t) \in \IR^m$ is the system input, and $y(t) \in \IR^m$ is the regulated system output. $A_{p} \in \mathcal{A} \subset \IR^{n \times n}$, $B_p \in \mathcal{B} \subset \IR^{n \times m}$, and $C_p \in \mathcal{C} \subset \IR^{m \times n}$ are the system matrices. $\omega_p \in \Omega \subset \IR^{m \times m}$, $\theta_p(t) \in \Theta \subset \IR^{n \times m}$, and $\sigma_p(t) \in \Delta \subset \IR^m$ are unknown system parameters. 
	Given a switching signal $p$, we assume that the sequence of finite switching time is $t_0,t_1,\ldots,t_i,\ldots$ with $t_0=0$.
	
	\begin{assum}
		The sets $\mathcal{A}$, $\mathcal{B}$, $\mathcal{C}$, $\Omega$, $\Theta$, and $\Delta$ are compact, convex polytopes.
		Without loss of generality, the sets $\Theta$ and $\Delta$ are assumed to contain $0$.
		$\Omega$ is assumed to be diagonally dominant, and without loss of generality, it is assumed to contain $\II$. 
		Define 
		\begin{align*}
			D_\theta \triangleq \max_{\theta \in \Theta} \norm{\theta}, \;
			D_\sigma \triangleq \max_{\sigma \in \Delta} \norm{\sigma}, \;
			D_\omega \triangleq \max_{\omega \in \Omega} \abs{\trace(\omega - \II)}.
		\end{align*}
		It is further assumed that $\sigma_i(t)$ and $\theta_i(t)$ are continuous and have (unknown) bounded derivatives for all $i\in \mathcal{I}$, i.e. 
		\begin{align}
		 \abs{\dot{\sigma}_i(t)} \leq d_\sigma, \quad \abs{ \dot{\theta}_i(t)} \leq d_\theta,\quad 	\forall i \in \mathcal{I}.
		\end{align}
	\end{assum}
	
	\begin{assum} \label{as::refsysstab}
		The switching signal $p$ has a {\it dwell time}, $\tau_d>0$, i.e. the switching times $t_1, \, t_2, \, \ldots$ satisfy the inequality $t_{k+1}-t_k \geq \tau_d$ for all $k$. However, the results derived for dwell-time switching also hold for the more general case of {\it average-dwell-time switching}. For the details of (average) dwell-time switching, see \cite{Liberzon}.  
	\end{assum}
	To analyze the performance of the adaptive system that will be presented in Section~\ref{sec:l1-archi}, we define a (non-adaptive) reference system that contains {\it perfect} knowledge of the parameters:
	\begin{equation}\label{eq:ref-sys}
		\begin{split}
			\dot{x}_\reft(t) &= A_p x_\reft(t) + B_p \big( \omega_p u_\reft(t) + \theta^\top_p(t) x_\reft(t) \\
			&\quad + \sigma_p(t) \big), \quad x_\reft(0) = x_0, \\
			u_\reft(s) &=  - \frac{D_0(s)}{s} \mu_\reft(s), 
		\end{split}
	\end{equation}
	where $\mu_\reft(t) \triangleq \omega_p u_\reft(t) + \theta^\top_p (t) x_\reft(t) + \sigma_p(t)-k_pr(t)$.
	The last equation in \eqref{eq:ref-sys} is equivalent to\footnote{We use the input-output mapping form instead of a transfer function form in \eqref{eq:uref-equi-filter-form} since the mapping is time-varying due to the existence of switching.} 
	\begin{equation}\label{eq:uref-equi-filter-form}
	 u_\reft = - \omega_p^{-1} \mcC_p \xi_\reft , 
	\end{equation}
	where $\xi_\reft(t) \triangleq \theta^\top_p (t) x_\reft(t) + \sigma_p(t) - k_pr(t)$ with $k_p$ being a feedforward gain for reference tracking, and the (time-invariant) mapping $\mcC_i$ has the transfer function form of
	\begin{equation} \label{eq:low-pass-filter}
		\mcC_i(s) \triangleq \omega_i (s\II+D_0(s)\omega_i)^{-1}D_0(s),
	\end{equation}
	which denotes a low-pass filter with the DC gain equal to an identity matrix, i.e. $\mcC_i(0) = \II$.
	From the first equation in \eqref{eq:ref-sys} and \eqref{eq:uref-equi-filter-form}, one can see that the reference control input tries to cancel the uncertainties within the bandwidth of the filter $\mcC_p$.
	This reference system provides the target performance of the \lone adaptive controller.
	
	Letting $(A_f,B_f,C_f,D_f)$ be a minimal realization of $D_0(s)$ with $n_f$ states, the reference system dynamics can be written in state-space form:
	\begin{align} 
		\label{eq::ch4refsysdyn2}
		\begin{bmatrix} \dot{x}_\reft(t) \\ \dot{x}_{f_1}(t) \\ \dot{x}_{I_1}(t) \end{bmatrix} &= \underbrace{\left[ \renewcommand\arraystretch{1.3} \begin{array}{c|c c}A_{p}+B_p\theta^\top_p & 0 & -B_p \omega_p \\ \hline B_f \theta^\top_p & A_f & B_f \omega_p \\ D_f \theta^\top_p & C_f & D_f \omega_p \end{array} \right]}_{\triangleq\bar{A}_p}  \begin{bmatrix} x_\reft(t) \\ x_{f_1}(t) \\ x_{I_1}(t) \end{bmatrix} \nonumber \\
		& \quad + \underbrace{\begin{bmatrix} B_p \\ B_f \\ D_f \end{bmatrix}}_{\triangleq\bar{B}_{p}} \sigma_p(t) - \underbrace{\begin{bmatrix} 0 \\ B_f \\ D_f \end{bmatrix}}_{\triangleq\bar{E}_p} k_p r(t) , \\
		\begin{bmatrix} u_\reft(t) \end{bmatrix} &= 
		\underbrace{\begin{bmatrix} 0 & 0 & -\II \end{bmatrix}}_{\triangleq\bar{C}} \begin{bmatrix} x_\reft(t) \\ x_{f_1}(t) \\ x_{I_1}(t) \end{bmatrix}, \quad \begin{bmatrix} x_\reft(0) \\ x_{f_1}(0) \\ x_{I_1}(0) \end{bmatrix} = \begin{bmatrix} x_0 \\ 0 \\ 0 \end{bmatrix}, \nonumber
	\end{align}
	or more compactly,
	\begin{align}
		\begin{split} \label{eq::ch4refsyscompact}
			\dot{\bar{x}}(t) &= \bar{A}_p \bar{x}(t) + \bar{B}_{p} \sigma(t) + \bar{E}_p r(t), \quad \bar{x}(0) = \bar{x}_{0}, \\
			u_\reft(t) &= \bar{C} \bar{x}(t),
		\end{split}
	\end{align}
	where $\bar{x}(t)\triangleq [x_\reft^\top(t), \; x_{f_1}^\top(t), \; x_{I}^\top(t)]^\top$, and $x_{f_1}(t)\in \mbR^{n_f}$ and $x_{I_1}(t)\in\mbR^m$ are the states of $D_0(s)$ and the integrator, respectively.
	
	
	\begin{assum} 
		The parameters of $D_0(s)$ are designed such that the reference system \eqref{eq::ch4refsysdyn2} (or \eqref{eq:ref-sys}) is stable for any switching signal with dwell time $\tau_d$. The following lemma provides a sufficient condition, which we assume to be satisfied.
	\end{assum}
	
	\begin{lem} \label{lem:dwelltime} If there exists symmetric matrices $ \bar{P}_i(\omega)$ such that for all $(\theta, \omega) \in \Theta \times \Omega$,
		\begin{equation}\label{eq:switching_stability_condition}
			\begin{split}
				\bar{P}_i&\geq \II, \ \forall i\in \mathcal{I}, \\
				\bar{A}_i^\top \bar{P}_i + \bar{P}_i \bar{A}_i  &\leq -\lambda \bar{P}_i, \ \forall i\in \mathcal{I}, \\
				\bar{P}_i &\leq \mu \bar{P}_j,\ 
				\forall i,j, \in \mathcal{I},
			\end{split}
		\end{equation}
		for some constants $\lambda>0$ and $\mu\geq 1$,
		then, for arbitrary $a^\star \in (0,1)$, the reference system \eqref{eq::ch4refsysdyn2} (or \eqref{eq:ref-sys}) is exponentially stable for any switching signal with dwell time 
		\begin{align} \label{eq::dwelltime}
			\tau_d \geq \frac{\ln(\mu)}{(1-a^\star) \lambda}.
		\end{align}
	\end{lem}
	
	The above lemma can be easily proved following \cite{Liberzon} using the switched Lyapunov function, defined as \begin{equation} \label{eq:Vp}
		V_p(\bar{x}(t)) \triangleq \bar{x}^\top(t) \bar{P}_p(\omega) \bar{x}(t).
	\end{equation}
	If we remove dwell time constraint for the switching signal, i.e. allowing arbitrary switching, a common Lyapunov function approach is usually employed to guarantee the stability. This is reflected in the following lemma. 
	\begin{lem} \label{lem::CLF}
		If there exists a constant symmetric matrix $\bar{P} > 0$ such that $\bar{A}_i^\top \bar{P} + \bar{P} \bar{A}_i \leq -\II$ for all $i \in \mathcal{I}$ and for all $(\theta, \omega) \in \Theta \times \Omega$, then the stability of the reference system \eqref{eq::ch4refsysdyn2} (or \eqref{eq:ref-sys}) is guaranteed by the Lyapunov function $V(\bar{x}(t)) \triangleq \bar{x}^\top (t)\bar{P}\bar{x}(t)$ under arbitrary switching. 
		
	\end{lem}
	\begin{rem}
		The arbitrary switching can be seen as a special case of dwell-time switching with $\tau_d =0, \mu = 1$.
		As the gain of $D_0(s)$, corresponding to the bandwidth of the low-pass filter $C(s)$ in \eqref{eq:low-pass-filter}, increases toward infinity, the effect of the uncertainty in the reference system vanishes. Thus, theoretically, $D_0(s)$ can always be designed to stabilize the reference system if the uncertainty-free reference system is stable. 
	\end{rem}
	

	In the remainder of this section, we present an adaptive control solution that ensures that the controlled system follows the switched reference system with quantifiable transient and steady-state performance bounds.
	
	\subsection{\lonew Adaptive Control Architecture}\label{sec:l1-archi}
	For the switched system in \eqref{eq::ch4sysdyn2}, we define the state predictor as
	\begin{align} 
		\begin{split}\label{eq:state_predictor}
			&\dot{\hat{x}}(t) = A_{p} \hat{x}(t) + B_p \left( \hat{\omega}(t) u(t) + \hat{\theta}^\top(t) x(t) + \hat{\sigma}(t) \right), \\
			& \hat{x}(t_i) = x(t_i),
		\end{split}
	\end{align}
	where $\hat{x}$ is the state of the predictor, and $\hat{\omega}$, $\hat{\theta}$, and $\hat{\sigma}$ are the parameter estimates governed by the adaptation laws:
	\begin{align}
		\begin{split} \label{eq:adaptive_law}
			\dot{\hat{\theta}}(t) &= \Gamma \Proj{\hat{\theta}(t)}{-x(t)\tilde{x}^\top(t)P_pB_p} , \quad \hat{\theta}(0) = 0, \\
			\dot{\hat{\sigma}}(t) &= \Gamma \Proj{\hat{\sigma}(t)}{-B_p^\top P_p \tilde{x}(t)} , \quad \qquad \hat{\sigma}(0) = 0, \\
			\dot{\hat{\omega}}(t) &= \Gamma \Proj{\hat{\omega}(t)}{-B_p^\top P_p \tilde{x}(t)u^\top(t)} , \quad \hat{\omega}(0) = \II,
		\end{split}
	\end{align}
	with $P_p$ being the (1,1) block of the $\bar{P}_p(\II)$ in \eqref{eq:Vp} (under the same partition as \eqref{eq::ch4refsysdyn2}), and $\tilde{x}(t) \triangleq \hat{x}(t) - x(t)$ being the prediction error. 
	
	The control law is defined as
	\begin{align} \label{eq:control_law}
		u(s) = -\frac{D_0(s)}{s} \mu(s),
	\end{align}
	where $\mu(s)$ is the Laplace transform of $\mu(t) = \hat{\omega}(t)u(t) + \hat{\theta}^\top(t)x(t) + \hat{\sigma}(t)-k_pr(t)$, respectively, $D_0(s)$ is a proper stable transfer function, $k_p$ is a reference scaling gain, and $r(t)$ is a bounded, piecewise continuous reference signal. For autonomous systems, the reference scaling gain is often the inverse of the $p$th system's DC gain, i.e. $k_p = -(C_p A_p^{-1}B_p )^{-1}$, but for human-controlled applications, this could be some other scaling or shaping of the human input. 
	Note that the control input will always be continuous due to the presence of the filter.

	\subsection{Analysis of the \lone Controller}
	Let $\tilde{\omega}_p(t) = \hat{\omega}(t) - \omega_p$, $\tilde{\theta}_p(t) = \hat{\theta}(t) - \theta_p(t)$, and $\tilde{\sigma}_p(t) = \hat{\sigma}(t) - \sigma_p(t)$. 
	Defining 
	\begin{equation}\label{eq:tileta_p-defn}
	  \tilde{\eta}_p(t) = \tilde{\omega}_p(t) u(t) + \tilde{\theta}_p^\top(t) x(t) + \tilde{\sigma}_p(t), 
	\end{equation}
	the prediction error dynamics can be formed from \eqref{eq::ch4sysdyn2} and \eqref{eq:state_predictor} as: 
	\begin{align} 
		\label{eq::ch4errdyn}
		\dot{\tilde{x}}(t) &= A_{p} \tilde{x}(t) + B_p \tilde{\eta}_p(t), \quad \tilde{x}(0) = 0.
	\end{align}
	
	\begin{lem} \label{lem::ch4prederr}
		The prediction error $\tilde{x}(t)$ is uniformly bounded,
		\begin{align} \label{eq:xtildebnd}
			\Linfnorm{\tilde{x}} \leq \sqrt{ \frac{ \beta_{\tilde{x}} }{\Gamma}}, 
		\end{align}
		for all $p \in \mathcal{P}$, where
		\begin{align*}
			\beta_{\tilde{x}} &\triangleq 4(D_\theta^2 + D_\sigma^2 + D_\omega^2) + 4 \lambda^{-1}(D_\theta d_\theta + D_\sigma d_\sigma).
		\end{align*}
	\end{lem}
	\begin{pf}
		Consider the Lyapunov function during the $i$-th time interval, $t \in [t_i, t_{i+1})$,
		\begin{align*}
			V_i(t) &= \tilde{x}^\top(t) P_i \tilde{x}(t) \\
			& \quad + \frac{1}{\Gamma} \left( \trace\left(\tilde{\theta}_i^\top(t) \tilde{\theta}_i(t) \right) + \trace\left(\tilde{\omega}_i^\top(t)\tilde{\omega}_i(t)\right) + \tilde{\sigma}_i^2(t) \right)
		\end{align*}
		Differentiating along trajectories of the prediction error dynamics \eqref{eq::ch4errdyn} and substituting the adaptive laws \eqref{eq:adaptive_law},
		\begin{align} \label{Vderiv}
			\dot{V}_i(t) &= \tilde{x}^\top(t) \left( A_{i}^\top P_i + P_i A_{i} \right) \tilde{x}(t) \nonumber\\
			&\quad + 2 \tilde{x}^\top(t) P_i B_i \left( \tilde{\theta}_i^\top(t) x(t) + \tilde{\omega}_i(t)u(t) + \tilde{\sigma}_i(t) \right) \nonumber \\
			&\quad + \frac{2}{\Gamma} \left( \tilde{\sigma}_i(t) \dot{\hat{\sigma}}(t) + \trace \left( \tilde{\theta}_i^\top(t) \dot{\hat{\theta}}(t) \right) + \trace \left( \tilde{\omega}_i(t)\dot{\hat{\omega}}(t) \right) \right) \nonumber \\
			&\quad - \frac{2}{\Gamma} \left( \tilde{\sigma}_i(t) \dot{\sigma}(t) + \trace \left( \tilde{\theta}_i^\top(t) \dot{\theta}(t) \right) \right), \nonumber \\
			&\leq - \lambda V_i(t) + \lambda \frac{\beta_{\tilde{x}}}{\Gamma}. 
		\end{align}
		Integrating the inequality \eqref{Vderiv} for all $t \in [t_i,t_{i+1})$, we have
		\begin{align} \label{Vint}
			V_i(t) \leq e^{-\lambda(t-t_i)}V_i(t_i) + \left( 1 - e^{-\lambda(t-t_i)} \right) \frac{\beta_{\tilde{x}}}{\Gamma}.
		\end{align}
		Since $\tilde{x}(t_i) = 0$, $
		V_i(t_i) \leq \left( D_\theta^2 + D_\omega^2 + D_\sigma^2 \right)/\Gamma,
		$
		which, along with \eqref{Vint}, implies that
		\begin{align*}
			\norm{\tilde{x}(t)}^2 \leq V_i(t) & \leq \frac{\beta_{\tilde{x}} }{\Gamma}.
		\end{align*}
		This bound holds for all $i$, and thus \eqref{eq:xtildebnd} holds.
		\hfill $\square$
	\end{pf}
	
	\begin{rem}
		(Effect of non-zero initialization error) In the above derivation, we have assumed that $\hat{x}(t_i) = x(t_i)$ (see \eqref{eq:state_predictor}), i.e. zero initialization error at each switching time, which leads to a uniform bound for the prediction error $\tilde{x}$, as shown in Lemma~\ref{lem::ch4prederr}. In the presence of non-zero initialization error, following the the steps outlined in \cite[Section~2.2.4]{L1Book}, one can show that that the performance bounds will have additional additive exponentially decaying terms. This proof is outside the scope of the current paper.
	\end{rem}
	
	Let $\xi_p(t) \triangleq\theta_p^\top(t)x(t) + \sigma_p(t)-k_pr(t)$. 
	Notice that the adaptive input in \eqref{eq:control_law} can be equivalently expressed as
	\begin{equation} \label{eq:control-law-filter-form}
   u = - \omega_p^{-1} \mcC_p (\xi_p+\tilde{\eta}_p), 
\end{equation}
where $\mcC_p$ is introduced in \eqref{eq:uref-equi-filter-form} and $\tilde \eta_p$ is defined in \eqref{eq:tileta_p-defn}	.

Applying \eqref{eq:control-law-filter-form} to the system dynamics in \eqref{eq::ch4sysdyn2}, we obtain
	\begin{align}
		\dot{\bar{x}}(t) &= \bar{A}_p \bar{x}(t) + \bar{B}_p \sigma_p(t) + \bar{E}_p \tilde{\eta}_p(t) - \bar{E}_p k_p r(t), 
	\end{align}
	where $\bar{x}(t) \triangleq [x^\top(t), x_{f_2}^\top(t), x_{I_2}^\top(t)]^\top$ with $x_{f_2}\in\mbR^{n_f}$ and $x_{I_2}\in \mbR^m$ being the states of $D_0(s)$ and the integrator, respectively, $\bar{x}(0) = [x_0^\top, \; 0, \; 0]^\top$, and $\bar A_p$, $\bar B_p$, $\bar E_p$ and $\bar C$ are defined in \cref{eq::ch4refsysdyn2}.
	The distance between the state of the reference system \eqref{eq::ch4refsysdyn2} and the actual system \eqref{eq::ch4sysdyn2}, $e(t) \triangleq x_\reft(t) - x(t)$, can be expressed as
	\begin{align} 
		\begin{split}\label{eq::ch4edyn}
			\bar{e}(t) &= \bar{A}_p  \bar{e}(t) + \bar{E}_p \tilde{\eta}_p(t), \quad \bar{e}(0) = 0,\\
			e_u(t) &= x_{I_1}(t),
		\end{split}
	\end{align}
	where $\bar{e}(t) \triangleq [e^\top(t), \; x_{f_1}^\top(t), \; x_{I_1}^\top(t)]^\top$, and $e_u(t) \triangleq u_\reft(t) - u(t)$. 
	From the prediction error dynamics \eqref{eq::ch4errdyn}, for all $i \in \mathcal{I}$, we have 
	\begin{align*}
		\tilde{\eta}_i(t) &= B_i^\dagger (\dot{\tilde{x}}(t) - A_i \tilde{x}(t)), \\
		D_0(s) \tilde{\eta}_i(s) &= s D_0(s) B_i^\dagger \tilde{x}(s) - D_0(s) B_i^\dagger A_i \tilde{x}(s),
	\end{align*}
	where $B_i^\dagger = (B_i^\top B_i)^{-1}B_i^\top$.
	Applying this to the error dynamics in \eqref{eq::ch4edyn},
	and letting 
	\begin{align*}
		\bar{H}_p & \triangleq \begin{bmatrix} -B_p C_f & -B_p D_f \omega_p \\ 0 & 0 \\ 0 & 0 \end{bmatrix}, & 
		\bar{J}_p \triangleq \begin{bmatrix} -D_f B_p^\dagger \\ -B_fB_p^\dagger A_p \\ -D_f B_p^\dagger A_p \end{bmatrix}, \\
		\bar{F}_p & \triangleq \begin{bmatrix} A_f & B_f\omega_p \\ C_f & D_f \omega_p \end{bmatrix}, 
		& \bar{G}_p \triangleq \begin{bmatrix} -B_fB_p^\dagger A_p \\ -D_f B_p^\dagger A_p \end{bmatrix}, 
	\end{align*}
	the reference error dynamics can be compactly written as
	\begin{align*}
		\begin{bmatrix} \dot{\bar{e}}(t) \\ \dot{\bar{x}}_{f_2}(t) \end{bmatrix} &= \begin{bmatrix} \bar{A}_p & \bar{H}_p \\ 0 & \bar{F}_p \end{bmatrix} \begin{bmatrix} \bar{e}(t) \\ \bar{x}_{f_2}(t) \end{bmatrix} + \begin{bmatrix} \bar{J}_p \\ \bar{G}_p \end{bmatrix} \tilde{x}(t), \\
		\begin{bmatrix} e_u(t) \end{bmatrix} &= \begin{bmatrix} \bar{C} & \bar{L}_p \end{bmatrix} \begin{bmatrix} \bar{e}(t) \\ \bar{x}_{f_2}(t) \end{bmatrix} -D_f B_p^\dagger \tilde{x}(t),
	\end{align*}
where $\bar{x}_{f_2}(t) = [x_{f_2}^\top(t), \; x_{I_2}^\top(t)]^\top$.
	
	\begin{thm} \label{thm::ch4ebnd}
		Consider the closed-loop adaptive system with the \lonew~controller defined via \cref{eq:state_predictor,eq:adaptive_law,eq:control_law} and the closed-loop reference system \eqref{eq::ch4refsysdyn2} (or \eqref{eq:ref-sys}). 
			Suppose that there exist $\bar{P}_i(\omega)$ ($i\in \mathcal{I}$) and some constants $\lambda>0$ and $\mu\geq 1$ such that the inequalities in \eqref{eq:switching_stability_condition} hold for all $(\theta, \omega) \in \Theta \times \Omega$, and the dwell time satisfies \eqref{eq::dwelltime}. Then, there exist positive constants $\kappa_1$ and $\kappa_2$ such that 
		\begin{align} \label{eq::ebnd}
			\Linfnorm{x_\reft-x} \leq \kappa_1 \frac{\beta_{\tilde{x}}}{\Gamma}, \quad
			\Linfnorm{u_\reft-u} \leq \kappa_2 \frac{\beta_{\tilde{x}}}{\Gamma}.
		\end{align}
	\end{thm}
	\begin{pf}
		Partitioning $\bar{P}_i$ of \eqref{eq:switching_stability_condition} along the same partition as $\bar{A}_i$ in \eqref{eq::ch4refsysdyn2}, we have
		\begin{align*}
			\bar{P}_i &= \left[ \begin{array}{c|c} P_i & R_i \\ \hline R_i^\top & S_i \end{array}\right].
		\end{align*}
		Further defining $Q_i \triangleq S_i - R_i^\top P_i^{-1} R_i$, and considering \eqref{eq:switching_stability_condition}, one obtains
		\begin{equation}
			\begin{split} \label{eq::switching_Q}
				Q_i&\geq \II, \ \forall i\in \mathcal{I}, \\
				\bar{F}_i^\top Q_i + Q_i \bar{F}_i  &\leq -\lambda Q_i, \ \forall i\in \mathcal{I}, \\
				Q_i &\leq \mu Q_j,\ 
				\forall i,j, \in \mathcal{I}.
			\end{split}
		\end{equation}

		Let $V_i(t) = \bar{e}^\top(t) \bar{P}_i \bar{e}(t) + \nu \bar{x}_{f_2}^\top(t) Q_{i} \bar{x}_{f_2}(t)$ on the time interval $[t_i,t_{i+1})$, where the scalar $\nu > 0$ satisfies
		\begin{align*}
			- \lambda a \bar{P}_i + \frac{1}{\nu \lambda a}\bar{P}_i \bar{H}_i Q_i^{-1} \bar{H}_i^\top \bar{P}_i < 0,
		\end{align*}
		with $a \in (0,a^\star)$. Such $\nu$ always exists since $\bar{P}_i > 0$.
		Differentiating $V_i(t)$ along the system trajectories, we have
		\begin{align}
			\dot{V}_i(t) &= \nu \bar{x}_{f_2}^\top(t) \left( \bar{F}_i^\top Q_i + Q_i \bar{F}_i \right) \bar{x}_{f_2}(t) \nonumber \\
			& \quad + \bar{e}^\top(t) \left( \bar{A}_i^\top \bar{P}_i + \bar{P}_i \bar{A}_I \right) \bar{e}(t) + 2 \bar{e}^\top(t) \bar{P}_i \bar{J}_i \tilde{x}(t) \nonumber \\
			& \quad + 2 \bar{e}^\top(t) \bar{P}_i \bar{H}_i \bar{x}_{f_2}(t) + 2 \nu \bar{x}_{f_2}^\top(t) Q_i \bar{G}_i \tilde{x}(t) \nonumber \\
			& \leq \begin{bmatrix} \bar{e} (t) \\ \bar{x}_{f_2}(t) \\ \tilde{x}(t) \end{bmatrix}^\top \begin{bmatrix} -\lambda \bar{P}_i & \bar{P}_i \bar{H}_i & \bar{P}_i \bar{J}_i \\ \bar{H}_i^\top \bar{P}_i & - \nu \lambda Q_i & \nu Q_i \bar{G}_i \\ \bar{J}_i^\top \bar{P}_i & \nu \bar{G}_i^\top Q_i & 0 \end{bmatrix} \begin{bmatrix} \bar{e}(t) \\ \bar{x}_{f_2}(t) \\ \tilde{x}(t) \end{bmatrix} \nonumber \\
			& \leq - (1-a) \lambda V_i(t) + g \norm{\tilde{x}(t)}^2, \quad \forall t\in [t_i,t_{i+1}) \label{eq::ch4Vdot}
		\end{align}
		where the last inequality follows from square completion, and the scalar $g$ is given by
		\begin{align*}
			g \triangleq \norm{ \begin{bmatrix} \bar{J}_i^\top \bar{P}_i & \nu \bar{G}_i^\top Q_i \end{bmatrix} \begin{bmatrix} -\lambda a \bar{P}_i & \bar{P}_i \bar{H}_i \\ \bar{H}_i^\top \bar{P}_i & -\nu \lambda a Q_i\end{bmatrix}^{-1} \begin{bmatrix}\bar{P}_i \bar{J}_i \\ \nu Q_i \bar{G}_i \end{bmatrix}}.
		\end{align*}
		Integrating the last inequality in \eqref{eq::ch4Vdot} and applying the bound on $\tilde{x}$ from \eqref{eq:xtildebnd}, we have
		\begin{align} 
			V_i(t) &\leq V_i(t_i) e^{-(1-a) \lambda (t-t_i)} + \int_{t_i}^t \!\! e^{-(1-a)\lambda (t-\tau)} g \norm{\tilde{x}(\tau)}^2 d\tau \nonumber \\
			& \leq V_i(t_i) e^{-(1-a) \lambda (t-t_i)} + \int_{t_i}^t \!\! e^{-(1-a)\lambda (t-\tau)} g \frac{\beta_{\tilde{x}}}{\Gamma} d\tau \nonumber \\
			& \leq V_i(t_i) e^{-(1-a) \lambda (t-t_i)} \nonumber \\
			& \quad + \frac{g}{(1-a)\lambda} \frac{\beta_{\tilde{x}} }{\Gamma} \left(1 - e^{-(1-a)\lambda(t-t_i)} \right) \nonumber \\
			& \leq \mu V_{i-1}(t_i) e^{-(1-a) \lambda (t-t_i)} \nonumber \\
			& \quad + \frac{g}{(1-a)\lambda} \frac{\beta_{\tilde{x}} }{\Gamma} \left(1 - e^{-(1-a)\lambda(t-t_i)} \right), \label{eq::ch4ISSlyap}
		\end{align}
	for any $t\in [t_i,t_{i+1})$,	where the last inequality follows from $V_{i}(t_i) \leq \mu V_{i-1}(t_i)$ implied by \eqref{eq:switching_stability_condition}, \eqref{eq::switching_Q}, and the definition of the Lyapunov functions.
		It follows that if at some switching time $t_i$,
		\begin{align} \label{eq::ch4lyapswitchbnd}
			V_{i-1}(t_i) & \leq \frac{g}{(1-a)\lambda} \frac{1-\mu^{-\frac{1-a}{1-a^\star}}}{1- \mu^{\frac{a-a^\star}{1-a^\star}}} \frac{\beta_{\tilde{x}}}{\Gamma} ,
		\end{align}
		then by applying the dwell time constraint \eqref{eq::dwelltime} and substituting \eqref{eq::ch4lyapswitchbnd} into \eqref{eq::ch4ISSlyap}, we have
		\begin{align*} 
			V_i(t_{i+1}) &\leq \frac{g}{(1-a)\lambda} \frac{1-\mu^{-\frac{1-a}{1-a^\star}}}{1- \mu^{\frac{a-a^\star}{1-a^\star}}} \frac{\beta_{\tilde{x}}}{\Gamma}.
		\end{align*}
		Thus if the condition in \eqref{eq::ch4lyapswitchbnd} holds for some $t_i$, it will hold for all $t_j$ with $j\geq i$. 
		Since $\bar{e}(0) = 0$ and $\bar{x}_{f_2}(0) = 0$, the inequality in \eqref{eq::ch4lyapswitchbnd} is satisfied at $t=0$, and thus at every switching time $t_i$ ($i=0,1,\cdots)$. 
		Applying the bound in \eqref{eq::ch4lyapswitchbnd} to \eqref{eq::ch4ISSlyap}, we obtain a uniform bound for the reference tracking error over the time interval $[t_{i},t_{i+1})$ as
		\begin{align}\label{eq::ebndderivation}
			\norm{ \begin{bmatrix} \bar{e}(t) \\ \sqrt{\nu} \bar{x}_{f_2} \end{bmatrix} }^2 \leq V_i(t) \leq \frac{\mu g}{(1-a)\lambda} \frac{1-\mu^{-\frac{1-a}{1-a^\star}}}{1- \mu^{\frac{a-a^\star}{1-a^\star}}} \frac{\beta_{\tilde{x}}}{\Gamma}.
		\end{align}		
		This holds for all $i \in \mathcal{I}$, and thus the bound is uniform for all $t$.
		We further notice that 
			\begin{align}
				e_u(t) &= \bar{C} \bar{e}(t) + \bar{L}_p \bar{x}_{f_2}(t) - D_f B_p^\dagger \tilde{x}(t). \label{eq:eubndderivation}
			\end{align} 
		From \eqref{eq::ebndderivation}, \eqref{eq:eubndderivation} and the bound on $\tilde{x}(t)$ in \eqref{eq:xtildebnd}, we can extract appropriate constants $\kappa_1$ and $\kappa_2$ such that \eqref{eq::ebnd} holds. The proof is complete. \hfill $\square$
	\end{pf}

	To summarize, if the reference system is stable (meets the conditions of Assumption \ref{as::refsysstab}), Theorem \ref{thm::ch4ebnd} guarantees that the states of the adaptive system follow those of the reference system with a bound proportional to $\Linfnorm{\tilde{x}}$. 
	From Lemma \ref{lem::ch4prederr}, the bound on $\Linfnorm{\tilde{x}}$ is proportional to the inverse of the square root of $\Gamma$. 
	Thus, by increasing the adaptation gain $\Gamma$, under the adaptive control law in \eqref{eq:state_predictor}, \eqref{eq:adaptive_law}, \eqref{eq:control_law}, the adaptive system can be made arbitrarily close to the reference system, which is stable.

	\section{Simulation Results}
	
	In the following example, the short period dynamics of a transport class aircraft\footnote{Based on NASA's Transport Class Model Aircraft Simulation \url{https://software.nasa.gov/software/LAR-18322-1}.} are considered during the approach phase. 
	As noted in Section \ref{sec::motivation}, pilots desire to sense changes in control effectiveness as the air speed varies, which motivates the usage of varying desired dynamics. 
	
	The flaps and gear are deployed at 162 knots, and the approach speed is taken to be 137 knots. 
	A different model is used at each 5 knot increment between these speeds.
	The system matrices are given in App. \ref{app::params}.
	As the vehicle slows down, the control effectiveness decreases as well, i.e. for the same control surface deflection, the generated moment is reduced. The states of the system are the angle of attack and the pitch rate.
	The input is the elevator deflection.
	
	In the simulation experiments, the plant models are switched every 20 seconds. 
	In order to excite the dynamics, a reference input of a 1 degree step is commanded to the elevator in the middle of each interval, with a return to a 0 degree reference command at the end of each interval.
	
	The design is performed using conservative bounds on the uncertainties, $\omega \in [0.5, 1.5]$, $\theta \in [-50,50] \times [-50,50]$ and $\sigma \in [-20, 20]$.
	These parameters correspond to a 50\% offset in the control effectiveness, a more than 50\% offset in state derivatives, and a trim offset of $20^\circ$.
	Aerodynamic parameters are typically known with higher accuracy than this. 
	Stability of the reference system was verified using a common Lyapunov function based on Lemma \ref{lem::CLF}, with $D_0(s) = 4\pi$. The adaptation gain was set to $\Gamma = 10^4$.
	
\begin{figure}
		\centering
		\includegraphics[clip, trim=0cm 0cm 0cm 0cm, width=7.4cm]{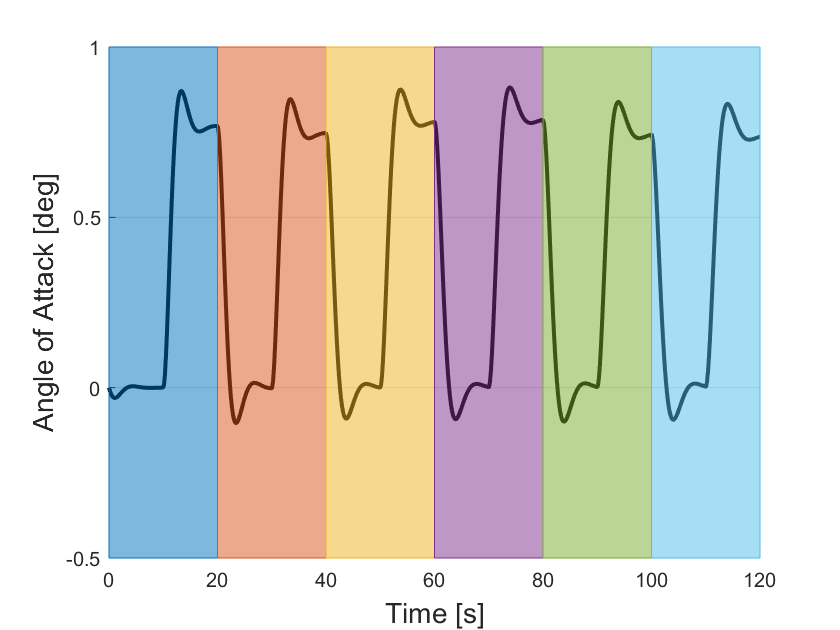}
		\vspace{-3mm}
		\caption{Angle of attack response with different design points}
		\label{fig::AoA}
	\end{figure}

\begin{figure}
		\centering
		\includegraphics[clip, trim=0cm 0cm 0cm 0cm, width=7.4cm]{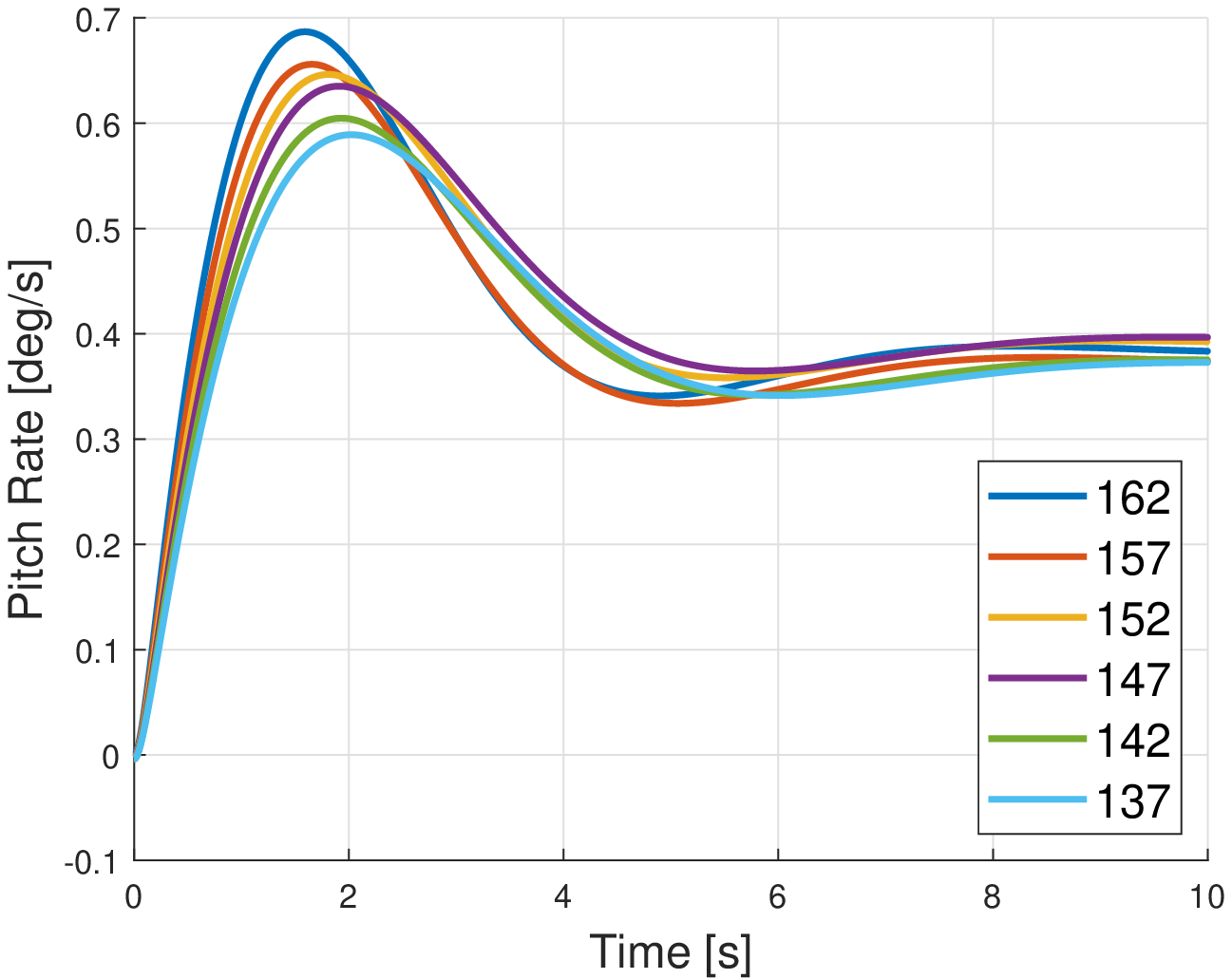}
		\includegraphics[clip, trim=0cm 0cm 0cm 0cm, width=7.4cm]{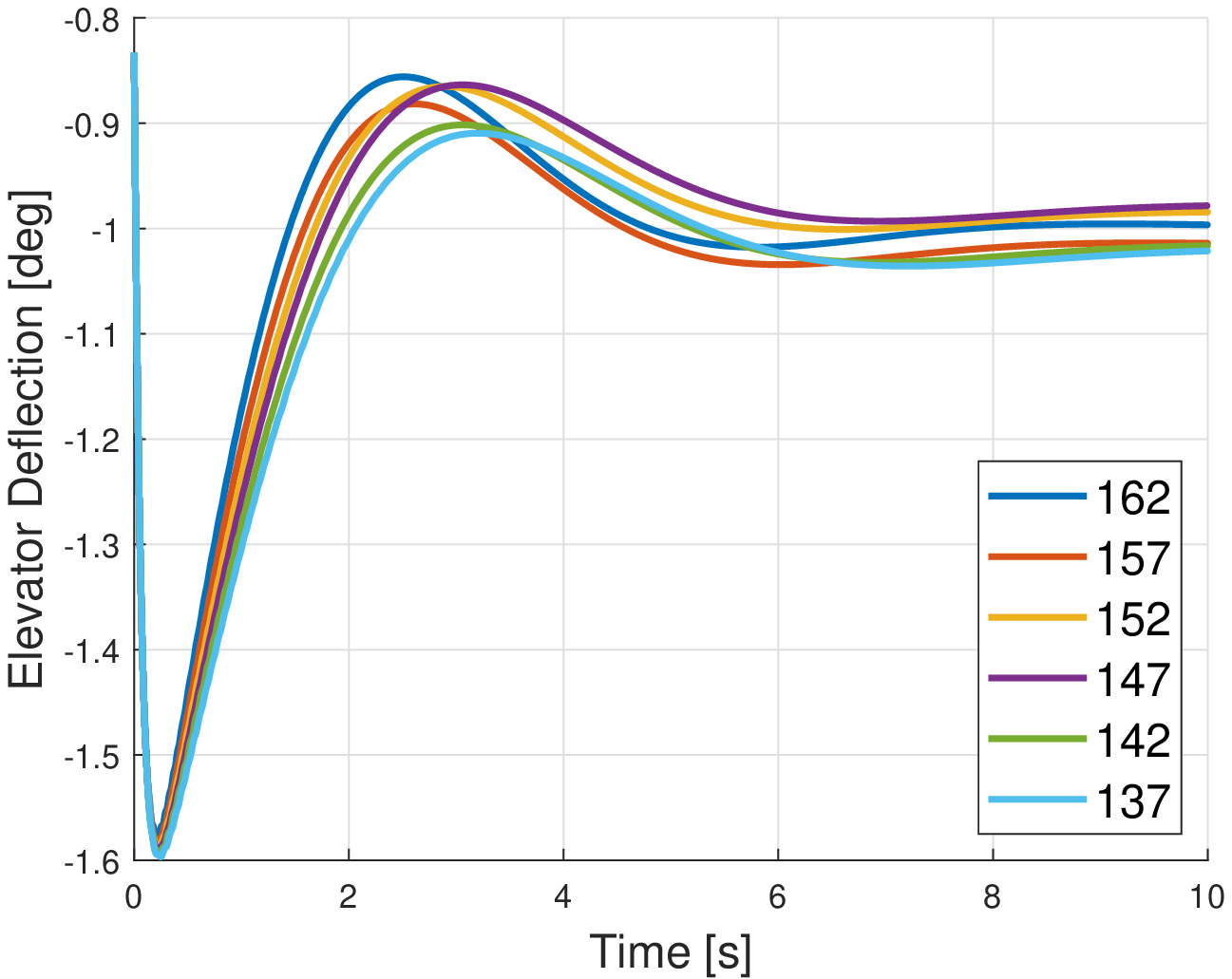}
		\vspace{-3mm}
		\caption{Comparison of pitch rates and elevator deflections at different velocities with different design points}
		\label{fig::pitchrate}
	\end{figure}

	The angle of attack response of both the adaptive system and the reference system is shown in Figure \ref{fig::AoA}, with uncertainty values of $\omega = 1.2$, $\theta = [-40; -40]$, and $\sigma = 1$. The response of the adaptive system is indistinguishable from that of the reference system. The pitch rate response during the upward rise in each interval is shown in Figure \ref{fig::pitchrate}, along with the elevator deflection. 
	The colored boxes in Figure \ref{fig::AoA} denote which system is active, and they correspond with the line colors of the overlaid responses in Figure \ref{fig::pitchrate}.
	From Figure \ref{fig::pitchrate}, speed of response, determined via the slope, or pitch acceleration, varies by 40\% over this speed range, enough to be noticeable by a pilot. 
	
\begin{figure}
		\centering
		\includegraphics[clip, trim=0cm 0cm 0cm 0cm, width=7.4cm]{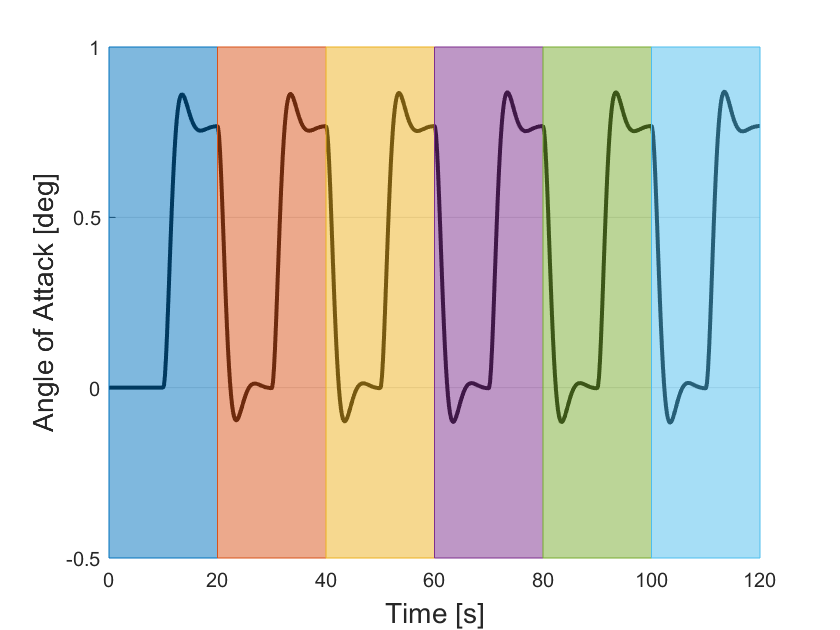}
		\vspace{-3mm}
		\caption{Angle of attack response for a single design point}
		\label{fig::AoA1}
	\end{figure}

\begin{figure}
		\centering
		\vspace{-8mm}
		\includegraphics[clip, trim=0cm 0cm 0cm 0cm, width=7.4cm]{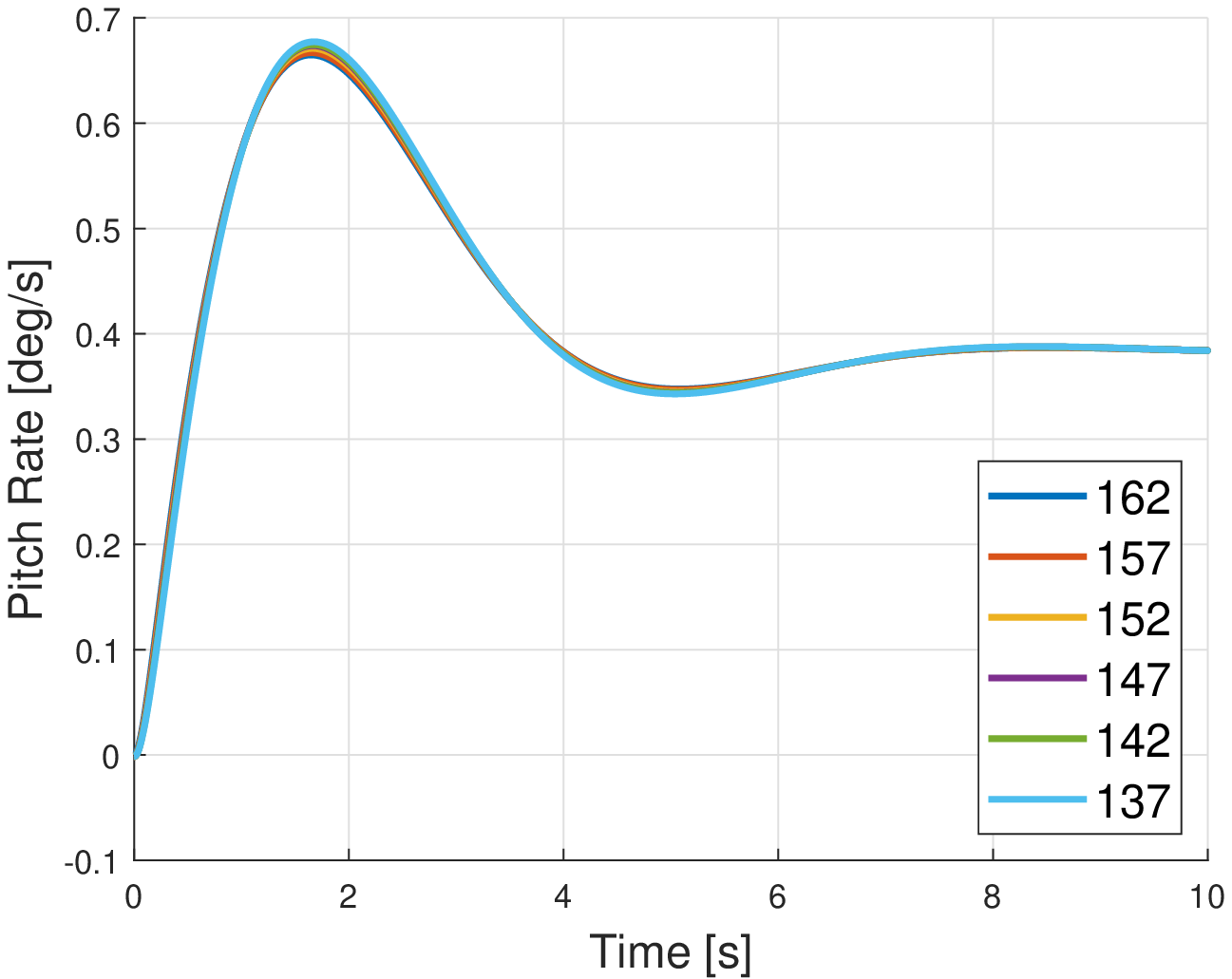}
		\includegraphics[clip, trim=0cm 0cm 0cm 0cm, width=7.4cm]{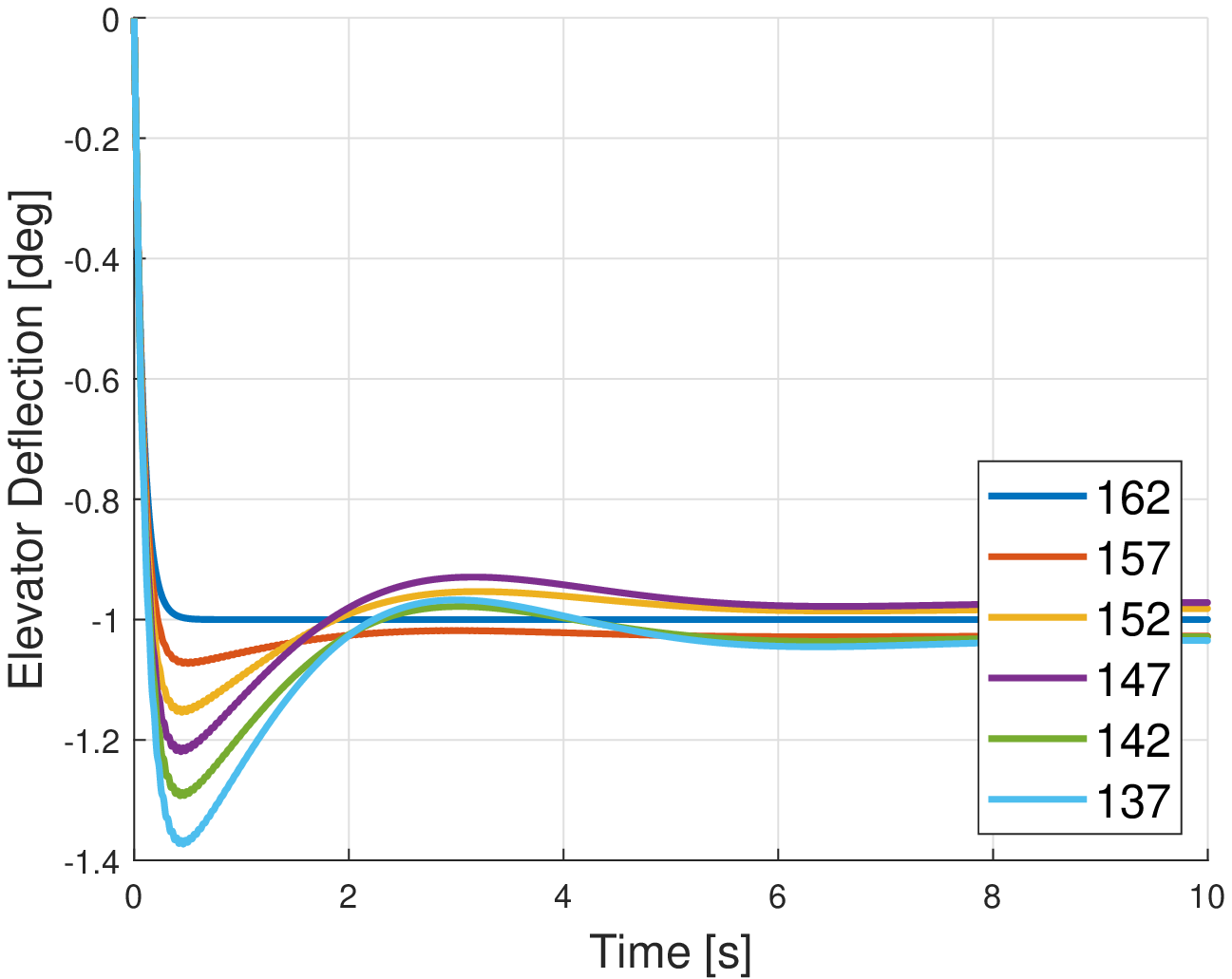}
		\vspace{-3mm}
		\caption{Comparison of pitch rates and elevator deflections at different velocities under a single design point}
		\label{fig::pitchrate1}
	\end{figure}

	For comparison, the same task was performed using a non-switched LTI reference system designed at the 162 knot condition. 
	For the fixed reference system, the only uncertainties were the differences between the fixed reference system and each of the plant's systems. 
	The angle of attack responses for this simulation are shown in Figure \ref{fig::AoA1}. 
	While Figure \ref{fig::AoA} shows some change in behavior as the reference model switches, Figure \ref{fig::AoA1} shows a consistent response, consistent with the theory in \cite{L1Book}. 
	The pitch rate response and elevator deflections during the upward rises are once again overlaid in Figure \ref{fig::pitchrate1} for comparison.
	The consistency of the response for a fixed reference system is clearly evident here.
	The expected variation in control effectiveness would not be apparent to pilots, requiring increased workload to maintain situational awareness. 
	

	\section{Conclusion}
	A framework for explicitly handling changes in the desired system dynamics during the design of an \lone adaptive controller has been presented. 
	This work was motivated by the common practice within the aerospace community of gain scheduling control laws. 
	It provides a means to mathematically verify the stability and robustness of a design despite changes in the desired system response. The proposed method was validated by the simulation of the short period dynamics of a transport class aircraft during the approach phase.
	
	
	\bibliography{draftBib}

\begin{thebibliography}{14}
\providecommand{\natexlab}[1]{#1}
\providecommand{\url}[1]{\texttt{#1}}
\providecommand{\urlprefix}{URL }
\expandafter\ifx\csname urlstyle\endcsname\relax
  \providecommand{\doi}[1]{doi:\discretionary{}{}{}#1}\else
  \providecommand{\doi}{doi:\discretionary{}{}{}\begingroup
  \urlstyle{rm}\Url}\fi

\bibitem[{Ackerman et~al.(2019)Ackerman, Puig-Navarro, Hovakimyan, Cotting,
  Duke, Carrera, McCaskey, Esposito, Peterson, and Tellefsen}]{Kasey2019}
Ackerman, K., Puig-Navarro, J., Hovakimyan, N., Cotting, M.C., Duke, D.J.,
  Carrera, M.J., McCaskey, N.C., Esposito, D., Peterson, J.M., and Tellefsen,
  J.R. (2019).
\newblock Recovery of desired flying characteristics with an {L1} adaptive
  control law: Flight test results on {Calspan's} {VSS} {Learjet}.
\newblock In \emph{AIAA SciTech Forum}. San Diego, CA.

\bibitem[{Ackerman et~al.(2016)Ackerman, Xargay, Choe, Hovakimyan, Cotting,
  Jeffrey, Blackstun, Fulkerson, Lau, and Stephens}]{Kasey2016}
Ackerman, K., Xargay, E., Choe, R., Hovakimyan, N., Cotting, M.C., Jeffrey,
  R.B., Blackstun, M.P., Fulkerson, T.P., Lau, T.R., and Stephens, S.S. (2016).
\newblock {L1} stability augmentation system for {Calspan's} variable-stability
  {Learjet}.
\newblock In \emph{AIAA SciTech Forum}. San Diego, CA.

\bibitem[{Davidson et~al.(1998)Davidson, Murphy, Lallman, Hoffler, and
  Bacon}]{harv}
Davidson, J.B., Murphy, P.C., Lallman, F.J., Hoffler, K.D., and Bacon, B.J.
  (1998).
\newblock High-alpha research vehicle lateral-directional control law
  description, analyses, and simulation results.
\newblock Technical Report {TP-1993-208465}, {NASA}.

\bibitem[{Gangsaas et~al.(2008)Gangsaas, Hodgkinson, Harden, Saeed, and
  Chen}]{busjet}
Gangsaas, D., Hodgkinson, J., Harden, C., Saeed, N., and Chen, K. (2008).
\newblock Multidisciplinary control law design and flight test demonstration on
  a business jet.
\newblock In \emph{AIAA Guidance, Navigation, and Control Conference}.
  Honolulu, HI.

\bibitem[{Gregory et~al.(2010)Gregory, Xargay, Cao, and Hovakimyan}]{AirSTAR}
Gregory, I.M., Xargay, E., Cao, C., and Hovakimyan, N. (2010).
\newblock Flight test of an {L1} adaptive controller on the {NASA} {AirSTAR}
  flight test vehicle.
\newblock In \emph{AIAA Guidance, Navigation, and Control Conference}. Toronto,
  Canada.

\bibitem[{Harris and Stanford(2018)}]{F35}
Harris, J.J. and Stanford, J.R. (2018).
\newblock F-35 flight control law design, development and verification.
\newblock In \emph{AIAA Aviation Forum}. Atlanta, GA.

\bibitem[{Heim et~al.(2018)Heim, Viken, Brandon, and Croom}]{L2F}
Heim, E.H.D., Viken, E.M., Brandon, J.M., and Croom, M.A. (2018).
\newblock {NASA}'s learn-to-fly project overview.
\newblock In \emph{AIAA Atmospheric Flight Mechanics Conference}. Atlanta, GA.

\bibitem[{Hovakimyan and Cao(2010)}]{L1Book}
Hovakimyan, N. and Cao, C. (2010).
\newblock \emph{{L1} Adaptive Control Theory}.
\newblock Society for Industrial and Applied Mathematics, Philadelphia, PA.

\bibitem[{Kersting and Buss(2017)}]{Kersting2017}
Kersting, S. and Buss, M. (2017).
\newblock Direct and indirect model reference adaptive control for
  multivariable piecewise affine systems.
\newblock \emph{IEEE Transactions on Automatic Control}, 62(11), 5634--5649.

\bibitem[{Liberzon(2003)}]{Liberzon}
Liberzon, D. (2003).
\newblock \emph{Switching in Systems and Control}.
\newblock Birkhäuser, Boston, MA.

\bibitem[{Puig-Navarro et~al.(2019)Puig-Navarro, Ackerman, Hovakimyan, Cotting,
  Duke, Carrera, McCaskey, Esposito, Peterson, and Tellefsen}]{Learjet2019}
Puig-Navarro, J., Ackerman, K., Hovakimyan, N., Cotting, M.C., Duke, D.J.,
  Carrera, M.J., McCaskey, N.C., Esposito, D., Peterson, J.M., and Tellefsen,
  J.R. (2019).
\newblock An {L1} adaptive stability augmentation system designed for
  {MIL-HDBK-1787} {Level} 1 flying qualities.
\newblock In \emph{AIAA SciTech Forum}. San Diego, CA.

\bibitem[{Sang and Tao(2012)}]{Sang2012}
Sang, Q. and Tao, G. (2012).
\newblock Adaptive control of piecewise linear systems: the state tracking
  case.
\newblock \emph{IEEE Transactions on Automatic Control}, 57(2), 522--528.

\bibitem[{Snyder et~al.(2019)Snyder, Zhao, and
  Hovakimyan}]{snyder2019switchingl1}
Snyder, S., Zhao, P., and Hovakimyan, N. (2019).
\newblock L1 adaptive control for switching reference systems: {Application} to
  flight control.
\newblock \emph{IFAC-PapersOnLine}, 52(16), 718--723.

\bibitem[{Yuan et~al.(2016)Yuan, Schutter, and Baldi}]{Yuan2016}
Yuan, S., Schutter, B.D., and Baldi, S. (2016).
\newblock Adaptive asymptotic tracking control of uncertain time-driven
  switched linear systems.
\newblock \emph{IEEE Transactions on Automatic Control}, 62(11), 5802--5807.

\end{thebibliography}
	\appendix
	\section{Projection Operator} \label{app::proj}
	Consider the following smooth convex function:
	\begin{align*}
		f(\theta) \triangleq \frac{(\epsilon_\theta+1)\theta^\top\theta - \theta_{max}^2}{\epsilon_\theta \theta_{max}^2},
	\end{align*}
	with $\theta_{max}$ being the norm bound imposed on the vector $\theta$, and $\epsilon_\theta > 0$ being a free parameter determining the projection tolerance. 
	The projection operator is defined as
	\begin{align*}
		\Proj{\theta}{y} \triangleq \begin{cases} y - \frac{\theta \theta^\top y f(\theta) }{\| \theta \|^2}  & \mbox{if } f(\theta) \geq 0 \mbox{ and } \theta^\top y > 0 \\ y & \mbox{otherwise}. \end{cases}
	\end{align*}
	For more details on the projection operator and its properties, see \cite{L1Book}.
	
	\section{Simulation Parameters} \label{app::params}
	\begin{align*}
		A_{162} &= \begin{bmatrix} -0.5301 & 0.9273 \\ -0.9106 & -0.6871 \end{bmatrix} & B_{162} &= \begin{bmatrix} -0.0009 \\ -0.0168 \end{bmatrix} \\
		A_{157} &= \begin{bmatrix} -0.5272 & 0.9289 \\ -0.8557 & -0.6580 \end{bmatrix} & B_{157} &= \begin{bmatrix} -0.0008 \\ -0.0154 \end{bmatrix} \\
		A_{152} &= \begin{bmatrix} -0.5201 & 0.9305\\ -0.7229 & -0.6279 \end{bmatrix} & B_{152} &= \begin{bmatrix} -0.0008 \\ -0.0141 \end{bmatrix} \\
		A_{147} &= \begin{bmatrix} -0.5168 & 0.9322 \\ -0.6618 & -0.5960 \end{bmatrix} & B_{147} &= \begin{bmatrix} -0.0007 \\ -0.0132 \end{bmatrix} \\
		A_{142} &= \begin{bmatrix} -0.5171 & 0.9339 \\ -0.6669 & -0.5637 \end{bmatrix} & B_{142} &= \begin{bmatrix} -0.0007 \\ -0.0123 \end{bmatrix} \\
		A_{137} &= \begin{bmatrix} -0.5147 & 0.9357 \\ -0.6219 & -0.5309 \end{bmatrix} & B_{137} &= \begin{bmatrix} -0.0006 \\ -0.0115 \end{bmatrix}
	\end{align*}
	
\end{document}